\documentclass[11pt]{article}
\usepackage{amsmath}
\usepackage{amsfonts}
\usepackage{amsthm}

\begin{document}

\def\N{\mathbb{N}}
\def\F{\mathbb{F}}
\def\Z{\mathbb{Z}}
\def\R{\mathbb{R}}
\def\Q{\mathbb{Q}}
\def\H{\mathcal{H}}

\parindent= 3.em \parskip=5pt

\centerline{\bf{A NOTE ON THE EQUALITY $\pi^2/6=\sum_{n=1}^{\infty} 1/n^2$}} 

\vskip 0.5 cm
\centerline{\bf{by A. Lasjaunias and J-P. Tran}}
\centerline {(Bordeaux, France)}
\vskip 0.5 cm
\hskip 1 cm \emph{In memory of James Stirling, great Scottish mathematician.} 
\vskip 0.5 cm 
\par In his famous treatise, entitled \emph{Methodus differentialis, sive tractatus de summatione et interpolatione serierum infinitarum} and dating from 1730, James Stirling (nicknamed The Venetian, 1692-1770) could prove the following :
$$ \text{Both series}\quad (I) \quad \sum_{n=1}^{\infty} 1/n^2 \quad \text{and}\quad (II) \quad \sum_{n=1}^{\infty} 3/[n^2\binom{2n}{n}] $$
$$\text{have the same limit}\quad l\in \R.$$
Eventhough, using the series $(II)$, he could compute a large number of digits for an approximation of this limit $l$, he apparently did not determine its exact value.(\emph{Here is a good example of a well known experience shared by many searchers, particularly in mathematics : giving up, being so close to a nice discovery.})
\newline A few years later, Leonhard Euler, studying the series $(I)$, could establish the equality $l=\pi^2/6$.
\par The only goal of this very short note is to show directly that the series $(II)$ does converge to $\pi^2/6$. 
\par We consider the function $g$ defined by $g(x)=Arcsinus(x)$ for $x \in ]-1;1[$. Recalling that $g$ is an odd function, from the equality $g'(x)=(1-x^2)^{-1/2}$ and by integration, we obtain the classical expansion
$$g(x)=\sum_{n=0}^{\infty} u_nx^{2n+1} \quad \text{where} \quad u_n=\frac{2^{-2n}}{2n+1}\binom{2n}{n}\quad \text{for}\quad  n\geq 0.\eqno{(1)}$$
We introduce the function $f$ defined by $f(x)=(g(x))^2$, for $x \in ]-1;1[$. We observe that $f$ is an even function with $f(0)=0$. Hence, using the previous expansion, we get
$$f(x)=\sum_{n=1}^{\infty}v_nx^{2n} \quad \text{where} \quad v_{n+1}=\sum_{i=0}^{i=n}u_iu_{n-i}\quad \text{for}\quad  n\geq 0.\eqno{(2)}$$ 
We shall compute the first and the second derivative of the function $f$. We have 
$$f'=2gg' \quad \text{and}\quad f''=2g'^2+2gg''.\eqno{(3)}$$
Since $g'(x)=(1-x^2)^{-1/2}$, we obtain $g''(x)=x(1-x^2)^{-3/2}$. Observing, from $(3)$, that $2g(x)=(1-x^2)^{1/2}f'$, we get the differential equation :
$$(1-x^2)f''(x)=2+xf'(x).\eqno{(4)}$$ 
Now we replace in $(4)$ the functions $f$, $f'$ and $f''$ by their corresponding expansions, following from $(2)$. Hence, identifying coefficients on both sides of $(4)$, we have 
$$v_1=1 \quad \text{and} \quad (2n+2)(2n+1)v_{n+1}=4n^2v_n \quad \text{for} \quad n\geq 1. \eqno{(5)}$$ 
We set $w_i=4i^2/(2i+2)(2i+1)$ for $i\geq 1$, hence we get  $v_{n+1}=\prod_{i=1}^{i=n}w_i$. Therefrom, we obtain $$v_{n+1}=2^{2n+1}(n!)^2/(2n+2)! \quad \text{for} \quad n\geq 0. \eqno{(6)}$$  
Since  $Arcsinus(1/2)=\pi/6$, we can write $\pi^2/36=f(1/2)= \sum_{n=1}^{\infty}v_n/2^{2n}$. Finally, introducing $(6)$, we have the desired result :
$$\pi^2/6= \sum_{n=1}^{\infty}3/[n^2\binom{2n}{n}].$$ 
\par A carefull reader will notice that no use has been made of the expansion for $g$ described in $(1)$. However, combining $(1)$,$(2)$ and $(6)$, we obtain, for $n\geq 0$, the following amazing identity :
$$(n+1)(2n+1)\binom{2n}{n}\sum_{i=0}^{i=n}\binom{2i}{i}\binom{2n-2i}{n-i}/[(2i+1)(2n-2i+1)]=2^{4n}.$$    
\vskip 1 cm
We are grateful to Jesus Guillera who, during a private conversation, pointed out the possibility of obtaining the limit of this second series by using inverse trigonometric functions.

\vskip 0.5 cm  
  
December 2023.

\end{document}